\documentclass[conference, 10 pt, conference]{ieeeconf}  

\IEEEoverridecommandlockouts                              
\usepackage{float}
\usepackage{amsmath}

\usepackage{amssymb}

\usepackage{amsthm}
\usepackage{algorithm} 
\usepackage{algpseudocode} 
\usepackage{xcolor}
\usepackage{bbm}     
\usepackage{graphicx} 
\usepackage[normalem]{ulem}

\newtheorem{lemma}{Lemma}

\newtheorem{definition}{Definition}
\newtheorem{problem}{Problem}
\newtheorem{assumption}{Assumption} 
\newtheorem{remark}{Remark} 
\newtheorem{example}{Example} 
\usepackage{dsfont}
\usepackage{booktabs}

\newcommand{\longthmtitle}[1]{\mbox{}{\textit{(#1):}}}

\title{\LARGE \bf Data-Enabled Predictive Control for Nonlinear Systems Based on a Koopman Bilinear Realization}

\author{Zuxun Xiong \, Zhenyi Yuan\,  Keyan Miao \, Han Wang \, Jorge Cort\'es \, Antonis Papachristodoulou
\thanks{Z. Xiong, K. Miao, and A. Papachristodoulou are with Department of Engineering Science, University of Oxford, Parks Road, Oxford, OX1 3PJ, U.K. {\tt\small \{zuxun.xiong, keyan.miao, antonis\}@eng.ox.ac.uk }}
\thanks{Z. Yuan and J. Cort\'es are with the Department of Mechanical and Aerospace Engineering, University of California, San Diego, La Jolla, CA 92093, USA, {\tt\small \{z7yuan, cortes\}@ucsd.edu}}
\thanks{H. Wang is with Department of Information Technology and Electrical Engineering, ETH Zurich, Zurich, Switzerland, {\tt\small hanwang@control.ee.ethz.ch}}
}

\begin{document}

\maketitle
\thispagestyle{empty}
\pagestyle{empty}

\begin{abstract}
This paper extends the Willems' Fundamental Lemma to nonlinear control-affine systems using the Koopman bilinear realization. 
This enables us to bypass the Extended Dynamic Mode Decomposition (EDMD)-based system identification step in conventional Koopman-based methods and design controllers for nonlinear systems directly from data. Leveraging this result, we develop a Data-Enabled Predictive Control (DeePC) framework for nonlinear systems with unknown dynamics. A case study demonstrates that our direct data-driven control method achieves improved optimality compared to conventional Koopman-based methods. Furthermore, in examples where an exact Koopman realization with a finite-dimensional lifting function set of the controlled nonlinear system does not exist, our method exhibits advanced robustness to finite Koopman approximation errors compared to existing methods. 


\end{abstract}

\section{Introduction}
The increasing complexity of modern dynamical systems have made accurate model identification increasingly challenging. This has motivated a growing interest in data-driven control methods, which aim to design controllers directly from data, bypassing explicit model identification \cite{RN12}.
Among these, direct methods based on Willems' Fundamental Lemma \cite{RN47} have received particular attention, as they enable control synthesis using sufficiently long collected trajectories without the need to perform explicit system identification. Following the comprehensive and seminal work that applied Willems' Fundamental Lemma to data-driven control \cite{RN10}, this line of research has developed fruitfully, leading to various extensions of traditional control methods.

For LTI systems, traditional controller design methods such as Linear Quadratic Regulator (LQR) and Model Predictive Control (MPC) were extended to data-driven versions, where explicit model identification is no longer required. Specifically, a data-driven formulation for LQR which is robust to noise is provided in \cite{RN45}. Different regularizers are also proposed to bridge direct and indirect methods, as well as certainty-equivalent and robust approaches \cite{RN181}. As for MPC, the Data-Enabled Predictive Control (DeePC) framework was proposed to design optimal controllers satisfying constraints directly on input-output data \cite{RN260}. Recursive feasibility and closed-loop stability are provided for this data-driven framework \cite{RN174}. Besides traditional control design methods, the direct data-driven method can also be applied to design more complex controllers like neural networks \cite{10886367}. 

Despite the remarkable success of direct data-driven control methods for linear systems, designing controllers for nonlinear systems remains an open question. A straightforward approach is to locally linearize the nonlinear system and then directly apply the data-driven methods developed for linear systems \cite{RN10}. Data collection method has been proposed to ensure that the data sufficiently captures the system behavior \cite{RN94}. 
Another line of work constructs dictionaries tailored to the system structure to lift nonlinear dynamics into a linear form \cite{RN71, berberich2020trajectory}, which requires prior knowledge of the system for accurate dictionary construction.
A recent work addresses system nonlinearity using regularized kernel methods within the DeePC framework to achieve optimal control by solving a non-convex optimization problem \cite{RN407}. An extension of the Fundamental Lemma to bilinear systems provides a theoretical foundation and suggests a pathway toward more general nonlinear control \cite{RN440}.
%

As an alternative to local or structure-dependent linearization, a promising direction is to exploit global linear representations of nonlinear dynamics. The Koopman operator offers such a perspective by lifting the original nonlinear system into a higher-dimensional space, where the dynamics evolve linearly \cite{koopman1931hamiltonian}. Although the Koopman operator was initially proposed for analyzing uncontrolled nonlinear systems, many existing applications to controlled systems incorporate the control input to the extended state space.
Koopman-based methods have been successfully applied in predictive control, demonstrating superior performance in terms of both prediction accuracy (compared to MPC based on other linearization techniques) and computational efficiency (compared to nonlinear MPC)~\cite{RN378}. They have also been used in robust control design with stability guarantee \cite{RN456}. We recommend \cite{RN436} to readers interested in a more comprehensive overview of Koopman-based control methods.
%
%
The Koopman operator enables the use of linear system theory and tools in the lifted space, making it a powerful framework for data-driven control of nonlinear systems.

The pioneer work~\cite{shang2024willems}, which serves as the main motivation for our treatment here, combines exact Koopman linear realizations with data-driven control, enabling direct control design for nonlinear systems bypassing the need to explicitly select lifting functions. 
 %
Recent studies, however, suggest that the Koopman bilinear realization, which represents the original nonlinear dynamics as bilinear rather than linear in the lifted space, offers significant advantages over the Koopman linear realization~\cite{RN450}. In particular, for control-affine systems, infinite-dimensional Koopman bilinear realizations can always be constructed, whereas a linear realization is not guaranteed to exist~\cite{RN441}. Motivated by this, we propose a direct data-driven control method based on Koopman bilinear realization for nonlinear system, in which the case of Koopman linear realization naturally emerges as a special case.

The first contribution of this paper is the Koopman bilinear realization-based formulation of a Fundamental Lemma for nonlinear systems. For nonlinear systems that admit exact
Koopman bilinear realization, we provide a data-based representation of any possible trajectory leveraging rich-enough input/state data. We also build connections of the derived result to existing works, showing that our result cover them as special cases. Secondly, we apply the proposed method to the direct data-driven predictive control problem for nonlinear systems. Compared to existing methods, we subsequently demonstrate the improved optimality of the proposed method, and particularly the enhanced robustness to approximation errors under inexact Koopman realizations.


\textbf{Notation:} We use $\mathbb{R}^{n\times m}$ to denote $n$-by-$m$ dimensional real matrices. Superscript $d$ is used to indicate a sequence corresponding to collected data, e.g., $u^d_{[1,T]}$ represents a collected input sequence from $u(1)$ to $u(T)$. $Y_{1,L,T-L+1}$ represents the $L$-order Hankel matrix of $T$-long signal $y$ starts from $y(1)$. We use blkdiag to represent block diagonal matrix. The Kronecker product $\otimes$ is defined for matrices $A\in\mathbb{R}^{m\times n}$ and $B\in \mathbb{R}^{p\times q}$ as $A \otimes B = 
\begin{bmatrix}
a_{11} B & \cdots & a_{1n} B \\
\vdots & \ddots & \vdots \\
a_{m1} B & \cdots & a_{mn} B
\end{bmatrix} \in \mathbb{R}^{mp\times nq}$. 

%
%


\section{Preliminaries and Problem Formulation}\label{sec:pre}

\subsection{Willems' Fundamental Lemma}
Consider a discrete-time linear time-invariant (LTI) system:
\begin{equation}\label{eq:linear system}
    x_{k+1} = Ax_k+Bu_k,
\end{equation}
where $x_k\in\mathbb{R}^{n_x}$ and $u_k\in\mathbb{R}^{n_u}$ are the system state and control input, respectively, and $A\in\mathbb{R}^{n_x \times n_x}$ and $B\in\mathbb{R}^{n_x\times n_u}$ are system matrices. The core idea of Willems' Fundamental Lemma is that any behavior of system~\eqref{eq:linear system} can be represented by a finite set of collected input/state data, and these data can be generated offline by a sequence of control input which is persistently exciting. To begin with, we first give the notion of $L$-persistently exciting data.

\begin{definition}\longthmtitle{$L$-Persistently exciting data for LTI systems}\label{def:PE-linear}
Consider the LTI system~\eqref{eq:linear system}, let $x^d_{[1,T]}$ be the system state sequence generated by a control input sequence $u_{[1,T]}^d$. If the matrix
\begin{equation*}
\left[\begin{array}{c}
X_{1, 1,  T-L+1} \\
U_{1, L, T-L+1}
\end{array}\right]
\end{equation*}
is full-row rank, then we say that the input/state data $(u_{[1,T]}^d,x^d_{[1,T]})$ is $L$-persistently exciting for system~\eqref{eq:linear system}.
\end{definition}

Note that, if we further ask the LTI system~\eqref{eq:linear system} to be controllable, then it suffices for the input signal $u^d_{[1,T]}$ alone to be $(L+n_x)$-persistently exciting \cite{RN47}, i.e., the Hankel matrix $U_{1,L+n_x,T-L-n_x+1}$ is full-row rank. In the next, we formally provide the Willems' Fundamental Lemma.
\begin{lemma}\longthmtitle{Willems' Fundamental Lemma~\cite{RN47,RN10}}\label{lem:WFL}
   Consider the LTI system~\eqref{eq:linear system}, let the input/state data $(u_{[1,T]}^d,x^d_{[1,T]})$ be $L$-persistently exciting for it, then
   \begin{itemize}
   \item[\textit{(i)}] There exists a $g \in \mathbb{R}^{T-L+1}$ such that any $L$-long input/state trajectory $(u_{[1,L]},x_{[1,L]})$ of the LTI system~\eqref{eq:linear system} can be expressed by
    \begin{equation}\label{eq:dd trajectory}
        \begin{bmatrix}
u_{[1,L]} \\
x_{[1,L]}
\end{bmatrix} = 
\begin{bmatrix}
U_{1,L,T-L+1} \\
X_{1,L,T-L+1}
\end{bmatrix} g.
    \end{equation}

\item[\textit{(ii)}] Conversely, for any $g \in \mathbb{R}^{T-L+1}$, 
\begin{equation}\label{eq:dd trajectory converse}
\begin{bmatrix}
U_{1,L,T-L+1} \\
X_{1,L,T-L+1}
\end{bmatrix} g
    \end{equation}
returns an $L$-long input/state trajectory of the LTI system~\eqref{eq:linear system}. 
\end{itemize}
\end{lemma}
In summary, Lemma~\ref{lem:WFL} allows one to characterize the behavior of the LTI system~\eqref{eq:linear system} through a finite set of persistently exciting data. This enables us to design controllers based on collected data without identifying the system.

\subsection{Koopman Bi/Linear Realization of Nonlinear Systems}
Consider a discrete-time nonlinear dynamical system
\begin{equation}\label{eq:nonlinear system}
    x_{k+1}=f(x_k,u_k),
\end{equation}
where $x_k\in\mathbb{R}^{n_x}$ and $u_k\in\mathbb{R}^{n_u}$ are the system state and control input, and $f: \mathbb{R}^{n_x \times n_u} \to \mathbb{R}^{n_x}$ is the nonlinear dynamic. For such a nonlinear system, the Koopman operator can be used to predict $x$ in the following form:
\begin{equation}\label{eq:KLR}
z_{k+1}=Az_k+Bu_k,\quad x_k=Cz_k,
\end{equation}
where $z\in\mathbb{R}^{n_z}$ is the lifted state defined as:
\begin{equation}
    z = \Psi(x):=
    \begin{bmatrix}\psi_1(x) & \cdots &
    \psi_{n_z}(x)
    \end{bmatrix}^{\top}.
\end{equation}
Here we refer the scalar functions $\psi_i(\cdot): \mathbb{R}^{n_x} \to \mathbb{R}$ for all $i=\{1,
\ldots,n_z\}$ as observables (or lifting/basis functions), and $\Psi(\cdot)$ the lifting function set. In general, the Koopman operator lifts the original nonlinear dynamic $f$ to a higher-dimensional space through observables, such that the lifted state $z$ follows a linear dynamic \cite{RN378}. Typically, $n_z\gg n_x$, and it could be infinity. We call \eqref{eq:KLR} a Koopman Linear Realization (KLR) of system \eqref{eq:nonlinear system}. If the lifted linear dynamics fully capture the behaviour of the original nonlinear system,i.e., there exists a bijective mapping between the original state and the lifted state under the Koopman realization,
we refer to it as an exact KLR.

However, even with an infinite-dimensional lifting function set, not every nonlinear system has an exact KLR. Many studies have considered using the Koopman operator to bilinearize nonlinear systems. This bilinearization approach has been demonstrated to provide a more accurate approximation of the system, particularly when the system dynamics are related to the terms involving the product of state and input terms, such as control affine systems. Specifically, for a nonlinear system \eqref{eq:nonlinear system}, we can find a lifting function set $\Psi(\cdot): \mathbb{R}^{n_x} \to \mathbb{R}^{n_z}$ that makes the lifted state $z=\Psi(x)$ follow a bilinear dynamic:
\begin{equation}\label{eq:KBR}
z_{k+1}=Az_k+Bu_k+H(z_k\otimes u_k),\quad x_k=Cz_k,
\end{equation}
where $H =\begin{bmatrix}
H_1 & \cdots & H_{n_z}
\end{bmatrix}^{\top}\in \mathbb{R}^{n_z \times n_un_z}$. Similarly, we call \eqref{eq:KBR} a Koopman Bilinear Realization (KBR) of system \eqref{eq:nonlinear system}. If the realization is exact, we call it an exact KBR. 


To get either a KLR or a KBR for a nonlinear system, one key step is applying the Extended Dynamic Mode Decomposition (EDMD) \cite{RN452} to identify the matrices $A$, $B$, $C$ and $H$ by solving a least squares problem based on collected trajectories of the original system. Once after the system identification, model-based controllers can be designed. We consider this approach indirect as far as control design is concerned, in that the model is identified first. Given that direct control design methods have demonstrated advantages over indirect approaches in linear systems, we seek to propose a control design framework for nonlinear system that bypasses EDMD and operates directly on data. In the following, we additionally provide a motivating example, showing that KBR outperforms KLR in general.
\begin{example}\longthmtitle{Prediction performance comparison}
Consider a control-affine system as
    \begin{equation}
        \begin{bmatrix}
        \dot{x}_1 \\ 
        \dot{x}_2
    \end{bmatrix}
    =
    \begin{bmatrix}
        0.3 x_1+u_1 \\ 
        0.2x_2-0.2x^2_1+\cos x_1 u_1+u_2
    \end{bmatrix},\label{eq:inexact KBR system}
    \end{equation}
which does not admit any finite-dimensional Koopman realization. We construct a finite-dimensional KLR (basis: monomials up to degree 4 and  $\begin{bmatrix}
\cos x_1 & \sin x_1
\end{bmatrix}^{\top}$) and a KBR (basis: $\begin{bmatrix}
x_1 & x_2 & \cos x_1 & \sin x_1
\end{bmatrix}^{\top}$) for it, respectively. We then compare the performance of two surrogate models with the true dynamics under the same input sequence after 20 and 400 time steps  (Figure \ref{fig:Inexact KBR prediction}). Although neither provides an exact realization of the system, the KBR demonstrates significantly better predictive performance with a smaller basis. While KLR has acceptable performance in the short term, its predictions gradually drift over time due to accumulated errors. This further enhances the motivation for using the more accurate KBR for data-driven control design.
\begin{figure}
  \centering
  \includegraphics[width=1\linewidth]{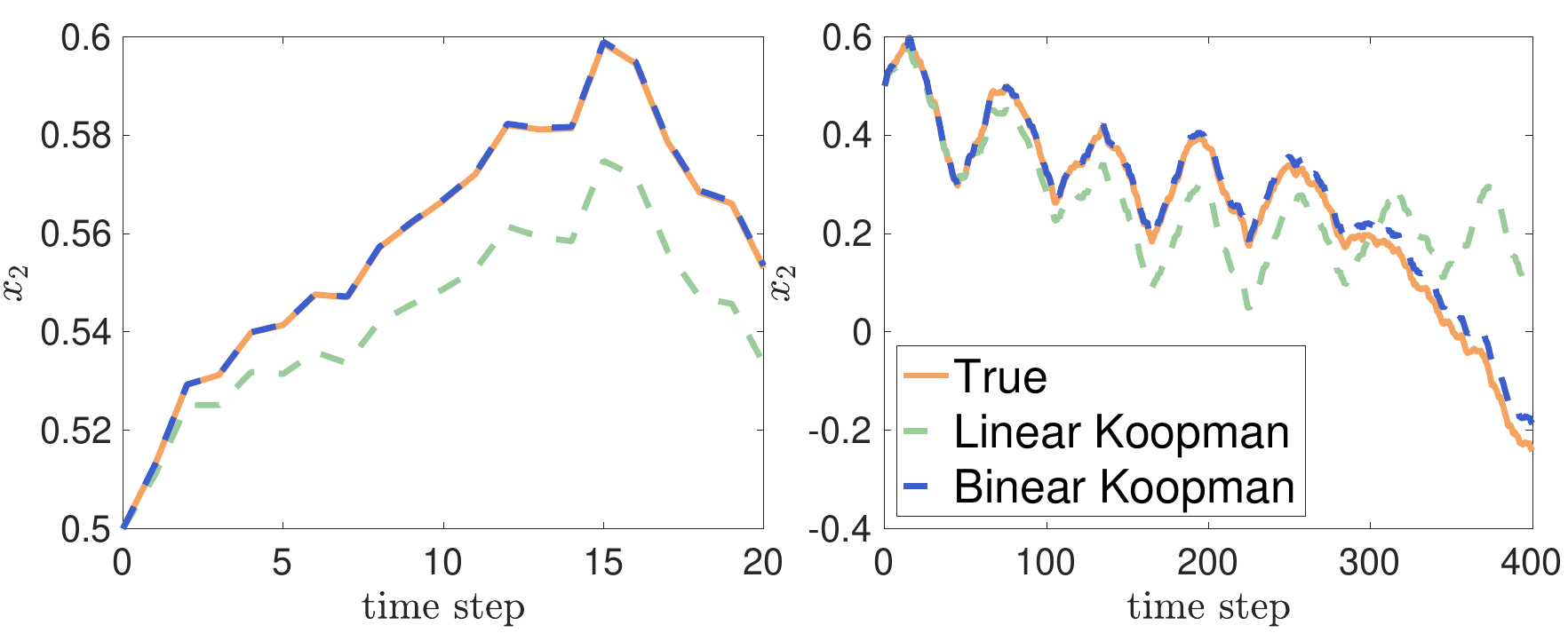}
  \caption{Prediction comparison between the KLR and the KBR (with different lifting function sets) of system \eqref{eq:inexact KBR system} after 20 prediction steps (left) and 400 prediction steps (right).}\label{fig:Inexact KBR prediction}
\end{figure}
\end{example} 

\subsection{Problem Formulation}
While direct data-driven control methods based on Willems' Fundamental Lemma have demonstrated notable success in linear systems, their extension to nonlinear systems remains limited.
This motivates us to solve:

\begin{problem}\label{prob:1}
    Extend Willems' Fundamental Lemma to nonlinear systems \eqref{eq:nonlinear system} via the Koopman operator, enabling direct data-driven control design for these systems.
\end{problem}

Moreover, designing predictive controllers for unknown nonlinear systems is an open problem: nonlinear optimization problems are inherently difficult to solve, and system identification for nonlinear dynamics can be inaccurate. To address these issues, we seek a direct data-driven solution that bypasses explicit model identification while enabling high-performance control. The second problem of our interest is:

\begin{problem}\label{prob:2}
Use the proposed Fundamental Lemma to solve following predictive control problem for nonlinear systems~\eqref{eq:nonlinear system} at time step $k$:
\begin{equation}\label{eq:nonlinear MPC}
\begin{aligned}
     \min_{u_i, \bar{x}_i} & \quad l_{N+1}(\bar{x}_{N+1}) + \sum_{i=1}^{N} \left[ l_i(\bar{x}_i) + u_i^\top R_i u_i + r_i^\top u_i \right] \\
     s.t. & \quad \bar{x}_{i+1} = f(\bar{x}_i, u_i), \quad i = 1, \ldots, N, \\
    & \quad c_{\bar{x}_i}(\bar{x}_i) + c_{u_i}^\top u_i \leq 1, \quad i = 1, \ldots, N, \\
    & \quad c_{\bar{x}_{N+1}}(\bar{x}_{N+1}) \leq 0, \\
    & \quad \bar{x}_1 = x_k.
\end{aligned}
\end{equation}
Here $\bar{x}$ represents the predicted state that only used within the MPC problem. It uses the true measured state $x_k$ of system~\eqref{eq:nonlinear system} at time step $k$ as the initial point. $N$ is the prediction horizon, $l(\cdot)$ is the nonlinear cost function on state $\bar{x}$, $R$ and $r$ are the quadratic and linear cost coefficients of the input, $c_{\bar{x}}(\cdot)$ and $c_u$ are constraints on state and input respectively.

\end{problem}

 
\section{Data-Driven Control Method for Nonlinear Systems Based on KBR}\label{sec:KDeePC}

\subsection{KBR-based Fundamental Lemma for Nonlinear System}
As emphasized earlier, the Koopman Bilinear Realization (KBR) offers higher accuracy than its linear counterpart. Building on this, we extend the Willems' Fundamental Lemma for bilinear systems proposed in \cite{RN440} based on the KBR setting, enabling its application to more general nonlinear systems. Theoretical results are provided based on following assumption:
\begin{assumption}\longthmtitle{Exact KBR of nonlinear system}\label{assump:exact finite KBR}
   There exists an exact KBR following dynamics \eqref{eq:KBR}, with $z=\Psi(x)\in\mathbb{R}^{n_z}, n_z<\infty$, for nonlinear system \eqref{eq:nonlinear system}.
\end{assumption}

Then we give a definition that characterizes when the collected input/state data is rich enough to represent any possible trajectory of nonlinear systems with exact KBR.

\begin{definition}\longthmtitle{$L$-persistently exciting data for nonlinear systems with exact KBR} \label{def:PE-bilinear}
Under Assumption \ref{assump:exact finite KBR}, let $x^d_{[1,T]}$ be the system state trajectory of \eqref{eq:nonlinear system} generated by a control input sequence $u_{[1,T]}^d$. Define lifted state trajectory and additional input sequence of its KBR~\eqref{eq:KBR} as $z^d_{[1,T]}:=\Psi(x^d_{[1,T]})$ and $v^d_{[1,T]}=(z^d \otimes u^d)_{[1,T]}$ respectively. If the matrix
\begin{equation}
\mathcal{G}_L(T) := \begin{bmatrix}
    Z_{1,1,T-L+1} \\
    U_{1,L,T-L+1} \\
    V_{1,L,T-L+1}
\end{bmatrix}
\in \mathbb{R}^{(n_z+n_uL+n_un_zL) \times (T-L+1)}
\end{equation}
is full-row rank, then we say that the input/state data $(u_{[1,T]}^d,x^d_{[1,T]})$ is $L$-persistently exciting for nonlinear system~\eqref{eq:nonlinear system} based on an exact KBR~\eqref{eq:KBR}.
\end{definition}

Now a fundamental lemma for nonlinear systems can be formulated based on Definition~\ref{def:PE-bilinear}:

\begin{lemma}\longthmtitle{KBR-based Fundamental Lemma for nonlinear systems}\label{lem:bilinear fundamental}
Under Assumption \ref{assump:exact finite KBR}, if the input/state data $(u_{[1,T]}^d,x^d_{[1,T]})$ collected from nonlinear system \eqref{eq:nonlinear system} satisfies the $L$-persistently exciting condition in Definition~\ref{def:PE-bilinear}, then the following statements hold.
\begin{itemize}
   \item[\textit{(i)}] Any input/state trajectory of system 
    \eqref{eq:nonlinear system} can be represented as
    \begin{equation}\label{eq:bilinear dd trajectory}
        \begin{bmatrix}
        x_{[1,L+1]} \\ 
        u_{[1,L]}
    \end{bmatrix}
    =
    \begin{bmatrix}
        X_{1,L+1,T-L+1} \\ 
        U_{1,L,T-L+1}
    \end{bmatrix}g
    \end{equation}
    for some $g\in \mathbb{R}^{T-L+1}$. 

\item[\textit{(ii)}]Conversely, any linear combination of data Hankel matrices, i.e., $\begin{bmatrix}
        X_{1,L+1,T-L+1} \\ 
        U_{1,L,T-L+1}
    \end{bmatrix}g$ is a valid input/state trajectory of system \eqref{eq:nonlinear system} if $g\in \mathbb{R}^{T-L+1}$ satisfies 
    \begin{equation}\label{eq:bilinear ddc cons}
        v_{[1,L]}=V_{1,L,T-L+1}g,
    \end{equation}
    where $v_i:=z_i\otimes u_i$ \footnote{Here $v,z,u$ are intermediate variables represent the additional input term, state, and input of KBR~\eqref{eq:KBR} for the nonlinear system \eqref{eq:nonlinear system} respectively. The whole trajectory $z_{[1,L+1]}$ and $u_{[1,L]}$ can be represented as $Z_{1,L+1,T-L+1}g$ and $U_{1,L,T-L+1}g$, respectively.}$, i=1,\ldots,L$.
\end{itemize}
\end{lemma}
\begin{proof}
Before proving this lemma, we introduce following matrices as in \cite{RN440}:
\[
\mathcal{O}_L := 
\begin{bmatrix}
    A \\
    A^2 \\
    \vdots \\
    A^L
\end{bmatrix}\in \mathbb{R}^{n_zL\times n_z},
\]

\[
\mathcal{P}_L = 
\begin{bmatrix}
    B & 0 & \cdots & 0 \\
    AB & B & \cdots & 0 \\
    \vdots & \vdots & \ddots & \vdots \\
    A^{L-1}B & A^{L-2}B & \cdots & B
\end{bmatrix}\in \mathbb{R}^{n_zL\times n_un_zL},
\]

\[
\mathcal{Q}_L = 
\begin{bmatrix}
    H & 0 & \cdots & 0 \\
    AH & H & \cdots & 0 \\
    \vdots & \vdots & \ddots & \vdots \\
    A^{L-1}H & A^{L-2}H & \cdots & H
\end{bmatrix}\in \mathbb{R}^{n_zL\times n_z}.
\]
Then we have $Z_{1,L,T-L+1}=[\mathcal{O}_L, \mathcal{P}_L, \mathcal{Q}_L] \mathcal{G}_L(T)$. 

    \textit{(i)} Since \( \mathcal{G}_L(T) \) is full-rank, we can find \( g\in \mathbb{R}^{T-L+1} \) satisfying:
\begin{equation}
    \begin{bmatrix}
    z(1) \\
    u_{[1, L]} \\
    v_{[1, L]}
\end{bmatrix}
= \mathcal{G}_L(T) g.\notag
\end{equation}

Based on the initial state $z(1)$ and input sequence $u_{[1,L]}$, we can represent the following state sequence as
\begin{equation}
\begin{aligned}
    z_{[2, L+1]} &= 
\begin{bmatrix} 
    \mathcal{O}_L & \mathcal{P}_L & \mathcal{Q}_L 
\end{bmatrix}
\begin{bmatrix}
    z(1) \\
    u_{[1, L]} \\
    v_{[1, L]}
\end{bmatrix} \\
&= 
\begin{bmatrix} 
     \mathcal{O}_L & \mathcal{P}_L & \mathcal{Q}_L
\end{bmatrix} 
\mathcal{G}_L(T) g = Z_{2, T-L+1} g.\notag
\end{aligned}
\end{equation}

Therefore any input/state trajectories of the KBR \eqref{eq:KBR} can be represented by:
\begin{equation}
\begin{bmatrix}
    z_{[1, L+1]} \\
    u_{[1, L]}
\end{bmatrix}
=
\begin{bmatrix}
    Z_{1, 1, T-L+1} \\
    Z_{2, L, T-L+1} \\
    U_{1, L, T-L+1}
\end{bmatrix} g
=
\begin{bmatrix}
    Z_{1, L+1, T-L+1} \\
    U_{1, L, T-L+1}
\end{bmatrix} g.
\notag
\end{equation}
Then the trajectories of nonlinear system~\eqref{eq:nonlinear system} can be represented by:
\begin{equation}
\begin{bmatrix}
    x_{[1, L+1]} \\
    u_{[1, L]}
\end{bmatrix}
\!=\!
\begin{bmatrix}
    \text{blkdiag}(C) & 0 \\
    0 & I
\end{bmatrix} 
\begin{bmatrix}
    z_{[1, L+1]} \\
    u_{[1, L]}
\end{bmatrix}
\!=\!
 \begin{bmatrix}
    X_{1, L+1, T-L+1} \\
    U_{1, L, T-L+1}
\end{bmatrix} g.
\notag
\end{equation}

\textit{(ii)} Consider the KBR~\eqref{eq:KBR} of system \eqref{eq:nonlinear system}: with $g$ satisfying \eqref{eq:bilinear ddc cons}, consider an initial point as 
\begin{equation}
z(1) = Z_{1, 1, T-L+1} g \notag,
\end{equation}
and an input sequence as
\begin{equation}
u_{[1, L]} = U_{1, L, T-L+1} g. \notag
\end{equation}

Then we have subsequent state as
\begin{equation}
\begin{aligned}
    z_{[2, L+1]} &= 
\begin{bmatrix} 
    O_L & P_L & Q_L 
\end{bmatrix}
\begin{bmatrix}
    z(1) \\
    u_{[1, L]} \\
    v_{[1, L]}
\end{bmatrix}\\
&=\begin{bmatrix} 
    O_L & P_L & Q_L 
\end{bmatrix} 
\mathcal{G}_L(T)g = Z_{2,L, T-L+1} g.
\end{aligned}\notag
\end{equation}
Therefore, $\begin{bmatrix}
    Z_{1, L+1, T-L+1} \\
    U_{1, L, T-L+1}
\end{bmatrix} g$ is a valid input/state trajectory of the KBR~\eqref{eq:KBR}. Subsequently, $\begin{bmatrix}
    X_{1, L+1, T-L+1} \\
    U_{1, L, T-L+1}
\end{bmatrix} g$ is a valid trajectory of \eqref{eq:nonlinear system}.
\end{proof}

\begin{remark}\longthmtitle{KLR-based formulation of Lemma~\ref{lem:bilinear fundamental}}\label{rem:lemma2 for linear}
    Consider a special case where an exact KBR of a nonlinear system~\eqref{eq:nonlinear system} exists with $H=0$ in \eqref{eq:KBR}, which means system~\eqref{eq:nonlinear system} admits an exact KLR. Then if the collected input/state data $u^d_{[1,T]}$ and $x^d_{[1,T]}$ is $L$-persistently exciting according to Definition~\ref{def:PE-linear}, we have: \textit{(i)} and \textit{(ii)} of Lemma \ref{lem:bilinear fundamental} still hold for system \eqref{eq:nonlinear system} without constraint \eqref{eq:bilinear ddc cons} being satisfied. This is because the bilinear term in the Koopman realization disappears. In this case, Lemma~\ref{lem:bilinear fundamental} reduces to the classical Fundamental Lemma (see Lemma~\ref{lem:WFL}), indicates that the knowledge of lifting function can be eliminated. In this sense, our result also recovers the result in~\cite{shang2024willems}, which also shows that for the data-based representation of trajectory of nonlinear systems which admit exact KLR, the need of lifting function can be bypassed.
\end{remark}

Before formulating a data-driven control design framework for a nonlinear system based on Lemma \ref{lem:bilinear fundamental} directly, it should be noticed that most nonlinear systems do not admit an exact finite-dimensional Koopman realization in practice. Only when the finite-dimensional space of observables is an invariant subspace of the Koopman operator can we exactly bi/linearizes the original dynamics \cite{RN436}. It is common to truncate the infinite observables set and use the finite part to provide an approximation. Specifically, we divide the infinite lifting function set into two parts as:
\begin{equation}
    \hat{z} = \begin{bmatrix}
        z\\
        z^\infty
    \end{bmatrix}=
    \begin{bmatrix}
        \Psi(x)\\
        \Psi^\infty(x)
    \end{bmatrix}=
    \begin{bmatrix}
    \begin{bmatrix}
        \psi_1(x)\\
    \vdots\\
    \psi_{n_z}(x)
    \end{bmatrix};
    \begin{bmatrix}
        \psi_{n_z+1}(x)\\
    \vdots\\
    \psi_{\infty}(x)
    \end{bmatrix}
    \end{bmatrix}.
\end{equation}
An exact KBR with infinite-dimensional $\hat{z}$ will be:
\begin{equation}
\begin{aligned}
    & \hat{z}_{k+1}=
    \begin{bmatrix}
        z_{k+1}\\
        z^{\infty}_{k+1}
    \end{bmatrix}=
    \begin{bmatrix}
        A_1 & A_2 \\
        A_3 & A_4
    \end{bmatrix}
    \begin{bmatrix}
        z_{k}\\
        z^{\infty}_{k}
    \end{bmatrix}+\begin{bmatrix}
        B_1 \\ B_2
    \end{bmatrix} u_k+\\
    &
    \begin{bmatrix}
        H_1 & H_2 \\
        H_3 & H_4
    \end{bmatrix}
    \begin{bmatrix}
        z_{k}\otimes u_k\\
        z^{\infty}_{k}\otimes u_k
    \end{bmatrix}, \quad x=\begin{bmatrix}
        C_1 & C_2
    \end{bmatrix}z.
\end{aligned}
\end{equation}
Therefore, we can still represent the dynamics of the finite-dimensional lifted state $z$ in a bilinear form by letting $A=A_1$, $B=B_1$, $C=C_1$, $H=H_1$, and $e_k^{\infty}=A_2z^{\infty}_{k}+H_2( z^{\infty}_{k}\otimes u_k)$:
\begin{equation}\label{eq:finite-KLR}
    z_{k+1}=Az_k+Bu_k+H(z_k\otimes u_k)+e_k^{\infty},\quad x_k=Cz_k.
\end{equation}
If the term in the lifting function set we discard is small enough or can be bounded, one can think of $e_k^{\infty}$ as a noise term. Some studies have provided finite error bounds for inexact Koopman bilinear realizations~\cite{RN453, RN460}.
%
%
We leave it as future work to formally consider it leveraging appropriate error metrics~\cite{MH-JC:24-scl}. Considering various studies have proposed a range of regularizers within the direct design setting to handle noise, model mismatch, and approximation errors \cite{RN260,RN266}, we propose a DeePC framework below with regularizers to effectively address the approximation errors inherent in most Koopman realizations.

\subsection{DeePC for Nonlinear System}
Using the Koopman operator, we can easily reformulate the nonlinear predictive control problem \eqref{eq:nonlinear MPC} into a linear one by choosing a special lifting function set $\Psi(x)$ that includes $\{x,l(x),c_x(x)\}$. Moreover, Lemma \ref{lem:bilinear fundamental} can be used to formulate a data representation for the system, thereby eliminating the need to identify an explicit model. The DeePC version of \eqref{eq:nonlinear MPC} can be formulated as:

\begin{subequations}\label{eq:DD KBR MPC}
\begin{align}
   &  \min_{z,u,v,g,\sigma_z} \quad 
    z^\top Qz + u^\top Ru+\lambda_z||\sigma_z||_2^2+\lambda_p p(g) \\
 &   s.t.    
     \begin{pmatrix}
    Z_{1,L+1,T-L+1} \\
    U_{1,L,T-L+1} \\
    \end{pmatrix}
    g = 
    \begin{pmatrix}
    \begin{bmatrix}
        z_{\text{ini}}\\
        z 
    \end{bmatrix}\\
    u 
    \end{pmatrix}+
    \begin{pmatrix}
   \begin{bmatrix}
        0 \\
        \sigma_z\\
    \end{bmatrix}\\
    0
    \end{pmatrix}, \label{eq:cons-ddc rep}\\
    & \quad v=V_{1,L,T-L+1}g, \label{eq:cons-bilinear1}\\
    & \quad v_i = z_i\otimes u_i, \quad i = 1, \ldots, N  \label{eq:cons-bilinear2}\\
    & \quad E_i z_i+F_i u_i \leq b_i, \quad i = 1, \ldots, N, \label{eq:cons-linear cons}\\
    & \quad z_{\text{ini}} = \Psi(x_k). \label{eq:cons-terminal}
\end{align}
\end{subequations}
Here $Q,q,R,r$ contain penalty coefficients on the inputs and states. Also, $\sigma_z$ is the slack variable to make the problem feasible even when there is no exact Koopman realization for the system, and $g\in \mathbb{R}^{T-N+1}$ is the vector we introduce for data-driven representation. The last two terms of the objective function are regularizers for the slack variable and $g$ respectively. The role and selection of different regularizers can be found in \cite{RN266}. Here we choose $p(g)=||g||^2_2$ to make the framework more robust to the Koopman approximation error. \eqref{eq:cons-ddc rep}-\eqref{eq:cons-bilinear2} is the data-driven representation of the original system and $U_d$, $Z_d$ and $V_d$ are Hankel matrices. Finally, $E,F,b$ are coefficients for constraints we derive from the original nonlinear MPC problem based on the chosen $\Psi(x)$, and $z_{\text{ini}}$ is the initial lifted state, which can be obtained from the measured current state $x_k$.

It is important to note that all data are collected from the original nonlinear system \eqref{eq:nonlinear system}. Specifically, we collect $T$-length input/state trajectories 
$x^d_{[1,T]}$ and $u^d_{[1,T]}$ respectively. We then compute $z^d_{[1,T]}$ based on the chosen basis $\Psi(\cdot)$, and subsequently construct $v^d_{[1,T]}$. Finally, the Hankel matrices are formed.

We note here that although Lemma~\ref{lem:bilinear fundamental} provides a data-based representation of input/state trajectory of original nonlinear system~\eqref{eq:nonlinear system}, we still formulate the DeePC problem in the input/state space of its KBR~\eqref{eq:KBR}. While it is equivalent for both when an exact KBR exists, we think that the formulation in~\eqref{eq:DD KBR MPC} is more robust when exact KBR does not exist, in the sense that the input/state space of the lifted bilinear system has higher dimension, and regularizers with higher dimension can potentially be more effective in handling such KBR errors. 
%
%
Also, the knowledge of lifting function set is better used in such formulation as we use the data of lifted state $Z$ instead of $X$.

\begin{remark}\longthmtitle{KLR-based DeePC}
    Following Remark \ref{rem:lemma2 for linear}, constraints \eqref{eq:cons-bilinear1} and \eqref{eq:cons-bilinear2} can be removed, resulting in a DeePC framework based on the KLR of the nonlinear system. Although the KLR is generally less accurate than the KBR in approximating nonlinear dynamics, it has demonstrated strong empirical performance in the control of many nonlinear systems. Besides, the KLR-based framework leads to a convex optimization problem while the KBR-based DeePC formulation includes bilinear constraints. This special case makes the proposed DeePC framework more flexible: when a finite-dimensional KLR provides a sufficiently accurate approximation of the nonlinear system, adopting this linear formulation enables more efficient computation and guarantees global optimality.
\end{remark}


\section{Numerical Examples}\label{sec:case}
We conduct three case studies to evaluate the effectiveness and advantages of the proposed method. For KLR-based approaches -- both the DeePC formulation and the EDMD-based MPC formulation -- the resulting optimization problems are convex with quadratic objective functions, which can be solved by OSQP solver \cite{osqp}. In contrast, the KBR-based formulations lead to non-convex problems due to bilinear constraints, and are solved using MATLAB's \texttt{fmincon}. All simulations are performed using MATLAB R2023b on a laptop with Apple M2 chip. 
\subsection{Nonlinear System with Inexact KLR}
A forced nonlinear Van der Pol oscillator \cite{RN362} is used to demonstrate the effectiveness of our method:
\begin{equation}
        \begin{bmatrix}
        \dot{x}_1 \\ 
        \dot{x}_2
    \end{bmatrix}
    =
    \begin{bmatrix}
        2x_2 \\ 
        -0.8x_1+2x_2-10x_1^2x_2+u
    \end{bmatrix}. \label{eq:VDP}
    \end{equation}
We start with the proposed KLR-based DeePC (KL-DeePC), and compare it with following methods:
    
    a) L-MPC: MPC based on linearization on current point;
    
    b) KL-MPC: MPC based on KLR with EDMD;
    
    c) KL-IO-DeePC: DeePC based on KLR with input-output data \cite{shang2024willems}.

Monomials up to degree 4 are selected as the basis functions for all methods. For other important parameters we have $N=12, Q=1, R=0.01, \lambda_g=0.01,\lambda_z=10$. 

\begin{figure}
  \centering
  \includegraphics[width=1\linewidth]{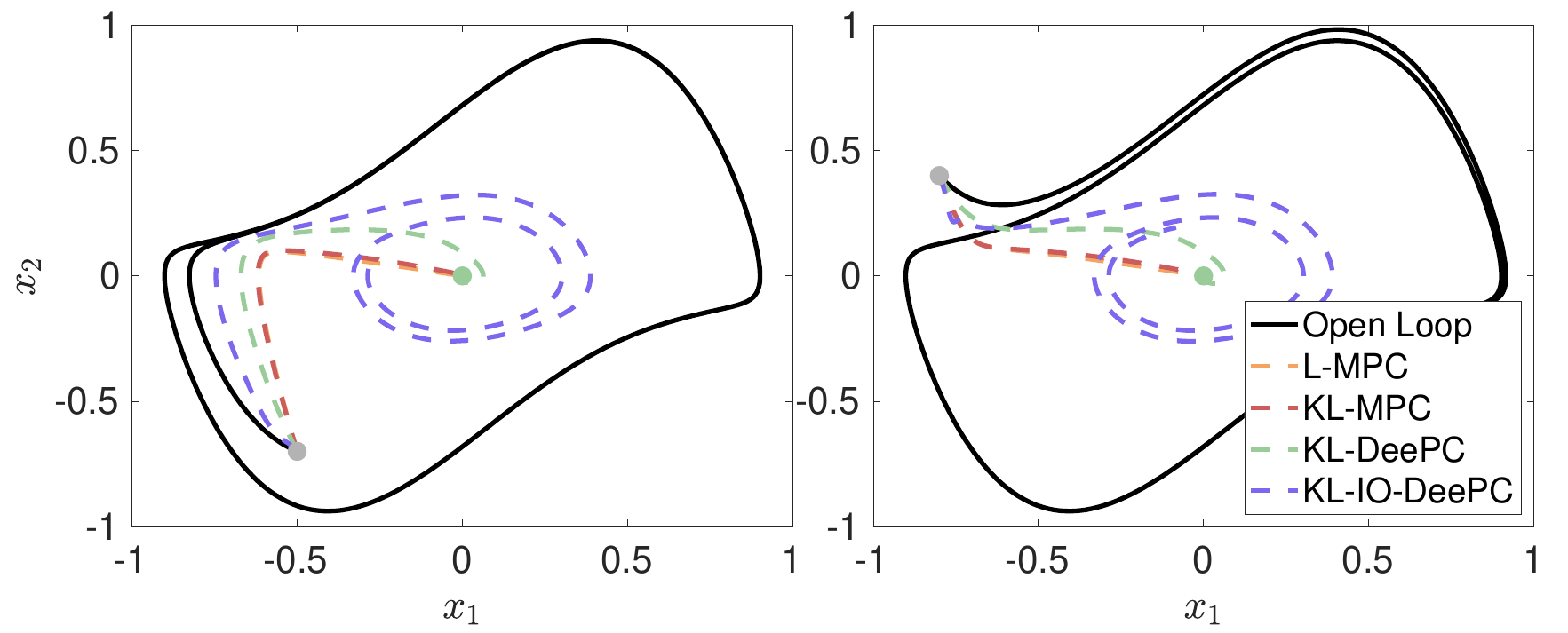}
  \caption{Open-loop and closed-loop optimal control trajectories of the Van der Pol oscillator under different control methods.}\label{fig:VDP traj}
\end{figure}

Figure \ref{fig:VDP traj} shows the trajectories from different initial conditions converging to the equilibrium point under different control methods over a simulation horizon of 1000 steps (0.02s for each step). The black curves represent the open-loop trajectories of the Van der Pol system. All methods successfully stabilize the system except KL-IO-DeePC, which heavily relies on the existence of an exact finite-dimensional Koopman realization. When such a realization cannot be obtained—which is the case in most practical scenarios—its performance degrades significantly. Notably, both the proposed method and KL-IO-DeePC use the same type of regularizer. As previously discussed, we attribute the improved performance of our approach to our DeePC formulation based on input/lifted state, which enables more effective use of the lifting function information. Furthermore, operating in a higher-dimensional lifted space allows the regularizer in our formulation to better mitigate the impact of Koopman approximation errors.

  \begin{table}
  \caption{Comparison of the total cost under different methods.}
    \centering
    \begin{tabular}{ccccc}
        \toprule
        $x(0)$ & L-MPC & KL-MPC & KL-IO-DeePC & KL-DeePC\\
        \midrule
        $[-0.5,-0.7]$ & 68.27 & 63.13 & 148.28 & \textbf{58.51}\\
        $[-0.8,0.4]$ & 81.07 & 75.13 & 121.17 & \textbf{47.73}\\
        \bottomrule
    \end{tabular}
    \label{tab:Inexact KLR cost}
\end{table}
Table~\ref{tab:Inexact KLR cost} summarizes the total costs under different control methods. The local linearization control method has a slightly higher cost than the Koopman linear MPC, as the Koopman operator enables global linearization of the system dynamics, leading to more accurate prediction and optimization of future states and inputs over the prediction horizon $N$. The proposed KL-DeePC achieves stabilization with the lowest cost, benefiting from its robustness to approximation errors.  By tuning the regularization penalty, the cost can be further optimized for different cases.

\begin{figure}
  \centering
  \includegraphics[width=0.8\linewidth]{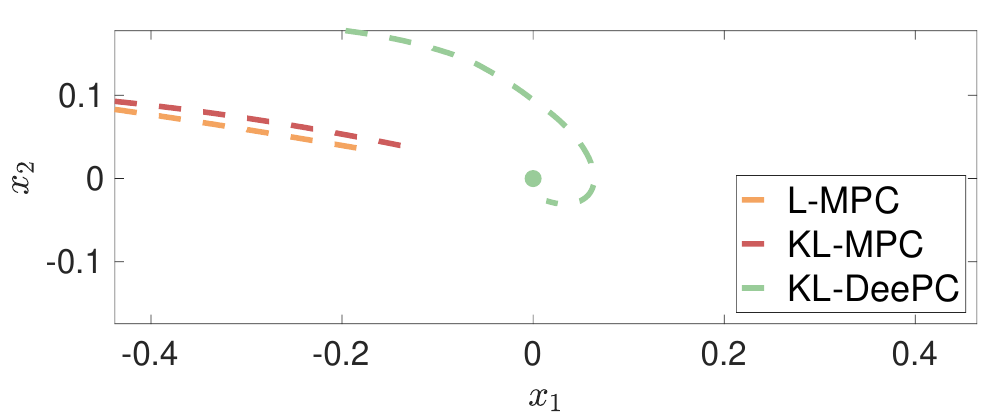}
  \caption{Closed-loop trajectories starting from $x(0)=[-0.8,0.4]$ under different methods after 400 steps (left) and 1000 steps (right).}\label{fig:conv1}
\end{figure}

Figure \ref{fig:conv1} shows the trajectories of system under different methods with initial state  $x(0)=[-0.8,0.4]$ after 400 steps. The KL-DeePC converges faster in this case.

\subsection{Nonlinear System with an Exact KBR}
 Consider a control-affine system from \cite{RN441}:
    \begin{equation}
        \begin{bmatrix}
        \dot{x}_1 \\ 
        \dot{x}_2
    \end{bmatrix}
    =
    \begin{bmatrix}
        0.3x_1+u_1 \\ 
        0.2x_2-0.2x^2_1+x_1^2u_1+u_2
    \end{bmatrix} \label{eq:exact KBR}
\end{equation}
There is no exact KLR for this system but we can construct an exact KBR using $\begin{bmatrix}
        x_1 & x_2 & x_1^2
    \end{bmatrix}^{\top}$ as the basis. We compare both direct and indirect design methods based on KLR and KBR respectively. 
    
    To better compare different methods, we consider both prediction horizon $N$ and control horizon $k_c$. As shown in Table \ref{tab:Exact KBR cost}, in the receding horizon strategy case: for methods based on KBR, as the bilinearization is exact, the proposed direct design method (KB-DeePC) shows the same performance as the one based on EDMD (KB-MPC); they both perform better in terms of the total cost compared to their counterparts based on KLR. We can obtain similar but more clear results from the second strategy. This is evident: as demonstrated in the motivating case, the approximation error of the KLR accumulates over longer horizons. Since system \eqref{eq:exact KBR} admits an exact KBR, the KB-DeePC and KB-MPC are equivalent. However, under the second control strategy, a slight difference is observed in the total cost between two methods. This discrepancy arises from the lack of global optimality guarantees in solving the non-convex optimization problem. Factors such as the choice of initial point may influence the solution, which highlights an important direction for future work.
    
\begin{table}
\caption{Comparison of the total cost under different methods.}
    \centering
    \begin{tabular}{ccccc}
        \toprule
        $(N,k_c)$ & KL-MPC & KL-DeePC & KB-MPC & KB-DeePC\\
        \midrule
        (10,1) & 39.07 & 39.00 & \textbf{38.78} & \textbf{38.78} \\
        (20,20) & 43.68 & 42.31 & \textbf{40.22} & 40.31 \\
        \bottomrule
    \end{tabular}
    \label{tab:Exact KBR cost}
\end{table}

\subsection{Nonlinear System with Inexact KBR}
Finally, we apply the proposed KB-DeePC to control system \eqref{eq:inexact KBR system}, where neither an exact KBR nor KLR exists, and compare its performance with both KB-MPC and KL-DeePC.
\begin{figure}
  \centering
  \includegraphics[width=0.95\linewidth]{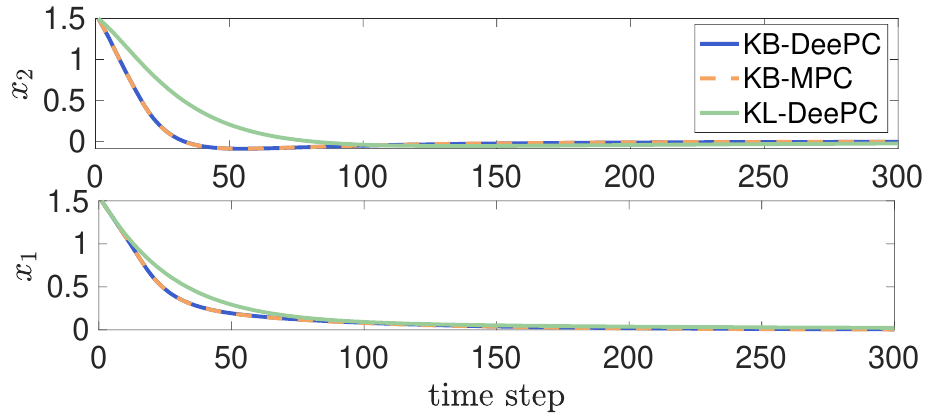}
  \caption{Comparison of state responses under different control methods.}\label{fig:Inexact KBR state}
\end{figure}
It can be seen from Figure \ref{fig:Inexact KBR state} that the proposed KB-DeePC achieves similar performance as the KB-MPC based on EDMD. This is because even if we cannot obtain an exact KBR for system \eqref{eq:inexact KBR system}, the bilinearization for this system based on a given basis is accurate enough, which is obvious from our motivating case in Section \ref{sec:pre}. 

In contrast, when comparing data-driven methods based on KBR and KLR respectively, the KBR-based approach clearly outperforms the KLR-based one. This is because that although both methods use the same Hankel matrix constructed from collected data, the KB-DeePC leverages the underlying bilinear dynamic, resulting in more accurate prediction and better control performance. 

\section{Conclusions And Future Work}\label{sec:conclusion}

We extend Willems’ Fundamental Lemma to nonlinear systems via Koopman bilinearization. For systems admitting an exact KBR, any input/state trajectory can be represented using PE data. The result also includes exact KLR systems as a special case, enhancing generality. To address the challenge that most nonlinear systems lack an exact finite-dimensional Koopman realization, we develop a DeePC framework where dynamics are represented using input/lifted state data. Combined with regularizers, this structure improves robustness to approximation errors. The proposed approach exhibits improved performance in a case study over traditional indirect and existing data-driven methods in nonlinear control. Future work includes formal error analysis and extensions to nonlinear networks.





\bibliographystyle{IEEEtran} 
\bibliography{main} 

\end{document}